\documentclass[11 pt]{amsart}

\usepackage{amscd, amssymb}

\newtheorem{pr}{Proposition}
\newtheorem{lm}{Lemma}
\newtheorem{tm}{Theorem}

\newcommand{\proj}{\mathbf P}

\newcommand{\barr}{\overline}
\newcommand{\rarr}{\rightarrow}
\newcommand{\oh}{{\mathcal{O}}}
\newcommand{\com}{\mathbb{C}}
\newcommand{\Q}{\mathbb{Q}}
\newcommand{\Z}{\mathbb{Z}}

\newcommand{\G}{\mathbf{G}}
\newcommand{\HH}{\Lambda}

\newcommand{\lan}{\langle}
\newcommand{\ran}{\rangle}
\newcommand{\eqq}{\stackrel{\sim}{=}}

\newcommand{\mgn}{\overline{M}_{g,n}}

\newcommand{\mgp}{\text{Map}_g}

\def\scup{\mathbin{\text{\scriptsize$\cup$}}}
\def\scap{\mathbin{\text{\scriptsize$\cap$}}}

\newcommand{\bpf}{\noindent {\em Proof.} }
\newcommand{\epf}{\qed \vspace{+10pt}}

\newcommand{\LL}{\mathbb{L}}
\newcommand{\hodge}{\mathbb{E}}

\newcommand{\ch}{\text{ch}}
\newcommand{\rk}{\text{rk}}
\newcommand{\CM}{\mathcal{M}}
\newcommand{\UM}{\mathcal{C}}

\begin{document}
\title{Hodge integrals and Gromov-Witten theory}
\author{C. Faber and R. Pandharipande}
\date{22 October 1998}
\maketitle

\pagestyle{plain}
\setcounter{section}{-1}
\section{\bf{Introduction}}
Let $\mgn$ be the nonsingular moduli stack of genus $g$, $n$-pointed,
Deligne-Mumford stable curves. For each marking $i$,
there is an associated cotangent line bundle $\LL_i\rarr \mgn$ with
fiber $T^*_{C,p_i}$ over the moduli point $[C, \ p_1, \ldots, p_n]$.
Let $\psi_i=c_1(\LL_i) \in H^*(\mgn, \Q)$. The integrals
of
products of the $\psi$ classes over $\mgn$ 
are determined by Witten's conjecture
(Kontsevich's theorem): their natural generating function
satisfies the Virasoro constraints [W], [K]. Let 
$\omega_C$ denote the dualizing sheaf of a curve $C$.
The Hodge bundle $\hodge\rarr \mgn$ is the rank $g$ vector bundle
with fiber $H^0(C, \omega_C)$ over
$[C, \ p_1, \ldots,p_n]$. Let $\lambda_j= c_j(\hodge).$
A Hodge integral over $\mgn$ is defined to be an integral
of products of the $\psi$ and $\lambda$
classes. It is the Hodge integrals that are studied here.

Hodge integrals arise naturally in
Gromov-Witten theory. There are two specific occurrences
which motivated this work.
First, let $X= \mathbf{G}/\mathbf{P}$ be
a compact algebraic homogeneous space. 
The virtual localization formula
established in [GrP] reduces all Gromov-Witten invariants (and
their descendents) of $X$ to  explicit graph sums involving only 
Hodge integrals over $\mgn$. For example, the classical Severi degrees
-- the numbers of degree $d$, genus $g$ algebraic plane
curves passing through $3d+g-1$ points -- are
Gromov-Witten invariants of $\proj^2$ and may be expressed
in terms of Hodge integrals. Formulas for 
Hodge integrals therefore play a role in
Gromov-Witten theory.

Second, let $X$ be an arbitrary nonsingular projective
variety of dimension $r$. Consider the stack $\overline{M}_{g,n}(X, 0)$
of stable {\em constant} maps from genus $g$, $n$-pointed curves to $X$.
There is a natural isomorphism:
\begin{equation}
\label{farst}
\overline{M}_{g,n}(X, 0) \eqq \mgn \times X.
\end{equation}
The virtual class $[\overline{M}_{g,n}(X, 0)]^{{vir}}$
is equal to $c_{rg}(\hodge^* \boxtimes T_X) \cap [ \overline{M}_{g,n}(X, 0)]$
via the identification (\ref{farst}). 
Hence, the degree $0$ Gromov-Witten invariants of $X$ involve only
the classical cohomology ring $H^*(X, \Q)$ and Hodge integrals
over $\mgn$. 
In [GeP], this observation is combined with the conjectural
Virasoro constraints 
of Eguchi, Hori, and Xiong [EHX] to yield 
simple formulas for certain Hodge
integrals. For example, the following relation is derived in
[GeP] as a consequence of the Virasoro constraints applied to
$\proj^1$:
\begin{equation} \label{lamg}
\int_{\mgn} \psi_1^{k_1} \dots \psi_n^{k_n} \lambda_g =
\binom{2g+n-3}{k_1,\ldots, k_n} b_g ,
\end{equation}
where $k_i\geq 0$ and
\begin{equation}
\label{qwq}
b_g = \begin{cases} 1 , & g=0 , \\ \displaystyle 
\int_{{\overline{M}}_{g,1}}
\psi_1^{2g-2} \lambda_g , & g>0.
\end{cases}
\end{equation}
The methods of [GeP] also yield conjectural relations for Hodge integrals
with a single $\lambda_{g-1}$ factor. The simplest of these predictions is:
for $g\geq 1$,
\begin{equation}
\label{ihop}
c_g = \int_{\overline{M}_{g,1}} \psi_1^{2g-1} \lambda_{g-1} = \Bigl(
\sum_{k=1}^{2g-1} \frac{1}{k} \Bigr) b_g - \frac12 
\sum_{\substack{g_1+g_2=g\\g_1,\,g_2>0}}
\frac{(2g_1-1)! (2g_2-1)!}{(2g-1)!} b_{g_1} b_{g_2}.
\end{equation}
Remarkably, the integrals $b_g$ seem to be unconstrained by the
degree 0 Virasoro conjecture.

More generally, it is natural to consider Hodge integrals
over stacks of stable maps $\mgn(X,\beta)$ for nonsingular
projective varieties $X$:
\begin{equation}
\label{hmint}
\int_{[\mgn(X,\beta)]^{vir}} 
\prod_{i=1}^n \psi_i^{a_i} \scup e_i^*(\gamma_i)
\scup \prod_{j=1}^g \lambda_j^{b_j}.
\end{equation} 
The classes $\psi_i$ here are the  cotangent line
classes on $\mgn(X,\beta)$, the maps $e_i$ are the evaluation
maps to $X$ corresponding to the markings, and the classes
$\gamma_i$ satisfy
$\gamma_i \in H^*(X, \Q)$. 
The gravitational descendents are the integrals 
(\ref{hmint}) for which all $b_j=0$ 
(no $\lambda$ classes appear). 
The first result proven in this paper is the
following Reconstruction Theorem.
\begin{tm}
\label{recon}
The set of Hodge integrals over moduli stacks of
maps to $X$ may be uniquely reconstructed from the
set of descendent integrals. 
\end{tm}
The method of proof is to utilize Mumford's
Grothendieck-Riemann-Roch calculations in [Mu]. Mumford's
results may be interpreted in Gromov-Witten theory to
yield differential equations for suitably
defined generating functions of Hodge integrals. 
A consequence of these equations is a direct
geometric construction of the $g=0$ relation $\tilde{L}_1$
which plays an important role in the proof
of the $g=0$ Virasoro constraints (see [EHX], [DZ], [Ge], [LiuT]).
As the required generating function
involves the Chern character of the Hodge bundle, it
seems quite difficult to obtain closed formulas for
the Hodge integrals (\ref{hmint}) via Theorem 1.
The reconstruction result was obtained 
in case $X$ is a point in [F2].

In order to find closed solutions in certain cases,
we introduce here a new method of obtaining relations among
Hodge integrals. The idea is to use the localization formula of
[GrP] in reverse: localization computations of known equivariant
integrals against $[\overline{M}_{g,n}(\mathbf{G}/\mathbf{P}, \beta)]^
{\text{vir}}$
yield relations among Hodge integrals over $\mgn$. 
A variant of this technique
is to compute an equivariant integral against the virtual class
via two different linearizations of the torus action. A relation
among Hodge integrals is then obtained by the two results of the
localization formula. A simpler case of these ideas provides motivation:
application of the Bott residue formula to integrals
over the Grassmannian 
yields nontrivial combinatorial identities when linearizations
are altered.

Hodge integrals over $\mgn$
also arise naturally in the study of tautological
degeneracy loci of the Hodge bundle. Formulas for these degeneracy loci
are used
here to find new relations among Hodge integrals. The geometry
involved is closely related to
classical curve theory: special linear series, Weierstrass points, and
hyperelliptic curves.

The main result of this paper is the following formula proven by 
the localization method together with a degeneracy calculation.
Define $F(t,k) \in \Q[k][[t]]$ by
$$F(t,k)= 1+ \sum_{g\geq 1} \sum_{i=0}^{g} t^{2g} k^i
\int_{\overline{M}_{g,1}} \psi_1^{2g-2+i} \lambda_{g-i}.$$
\begin{tm}
\label{mayne}
$$F(t,k)= \Big( \frac{t/2}{\sin(t/2)} \Big)^{k+1}.$$
\end{tm} 
\noindent
In particular, the integrals $b_g$ and $c_g$ are determined by:
\begin{equation}
\label{bgcg}
\sum_{g\geq 0} b_g t^{2g} = F(t,0)=  \Big( \frac{t/2}{\sin(t/2)} \Big),
\end{equation}
$$\sum_{g\geq 1} c_g t^{2g}= \frac{\partial F}{\partial k}(t,0) =
\Big( \frac{t/2}{\sin(t/2)} \Big) 
\cdot \log\Big( \frac{t/2}{\sin(t/2)} \Big).$$
D. Zagier has provided us with a proof of  
the Virasoro prediction
 (\ref{ihop}) from (\ref{bgcg}) and identities among
Bernoulli numbers. 
M. Shapiro and A. Vainshtein informed us of another approach
to Theorem 2 [ELSV], see also [SSV].

Theorem 2 has a direct application in Gromov-Witten theory to
a multiple cover formula for Calabi-Yau 3-folds. 
Under suitable conditions,
the integral
\begin{equation}
\label{exxx}
C(g,d)= \int _{[ \barr{M}_{g,0}(\proj^1,d)]^{\it{vir}}} 
c_{\rm{top}} (R^1 \pi_* \mu^* N)
\end{equation}
is the contribution to the
genus $g$ Gromov-Witten invariant of a Calabi-Yau 3-fold 
of multiple
covers of a fixed rational curve (with normal bundle
$N=\oh(-1)\oplus \oh(-1)$). 
The genus 0  case is determined by the Aspinwall-Morrison formula
$$C(0,d)=1/d^3,$$
[AM], [Ma], [V]. The genus 1 case was computed in physics [BCOV]
and mathematics [GrP] to yield $$C(1,d)= 1/12d.$$  Virtual localization
and Theorem 2 determine this multiple cover contribution in the general
case.

\begin{tm} For $g\geq 2$, $$C(g,d)= 
 \frac{|B_{2g}| \cdot d^{2g-3}}{2g\cdot (2g-2)!} =
|\chi(M_g)| \cdot \frac{d^{2g-3}}{(2g-3)!},$$
where $B_{2g}$ is the $2g^{th}$ Bernoulli number and
$\chi(M_g)= B_{2g}/2g(2g-2)$ is the Harer-Zagier formula for the orbifold
Euler characteristic of $M_g$.
\end{tm}
\noindent Theorem 3 was conjectured in [GrP] from data obtained
from the Hodge integral algorithm of [F2].

Another consequence of Theorem 2 is the determination
of the following Hodge integral.
\begin{tm} For $g\ge2$,
$$
\int_{\overline{M}_g}\lambda_{g-1}^3=\frac{|B_{2g}|}{2g}
\frac{|B_{2g-2}|}{2g-2}\frac1{(2g-2)!}.
$$
\end{tm}
\noindent The genus $g\geq 2$, degree 0 Gromov-Witten invariant of a
Calabi-Yau 3-fold $X$
is simply 
$$<1>_{g,0}^X= (-1)^g \frac{\chi}{2} \int_{\overline{M}_g} \lambda_{g-1}^3,$$
where $\chi$ is the topological Euler
characteristic of $X$ (see [GeP]).
Theorem 4 was conjectured previously in [F1],
and was recently derived in string theory
by physicists [MM], [GoV].
It implies Conjecture 1 in [F2].

We mention finally an interesting connection between
Gromov-Witten theory and the intrinsic geometry of $M_g$
via the Hodge integrals.
The ring $\mathcal{R}^*(M_g)$ of tautological Chow classes in
$M_g$ has been conjectured in [F1] to be a Gorenstein
ring with socle in degree $g-2$. The Hodge integrals 
\begin{equation}
\label{carcon}
\int_{\overline{M}_{g,n}} 
\psi_1^{k_1} \dots \psi_n^{k_n} \lambda_g \lambda_{g-1}
\end{equation}
determine the top intersection pairings in $\mathcal{R}^*(M_g)$. The study of
$\mathcal{R}^*(M_g)$ in [F1] led to a simple combinatorial conjecture for the
integrals (\ref{carcon}):  
\begin{equation} \label{lamgg}
\int_{\mgn} \psi_1^{k_1} \dots \psi_n^{k_n} \lambda_g \lambda_{g-1}
= \frac{(2g+n-3)! (2g-1)!!}{(2g-1)!\prod_{i=1}^n (2k_i-1)!!}
\int_{\overline{M}_{g,1}} \psi_1^{g-1} \lambda_g \lambda_{g-1},
\end{equation}
where $g\geq 2$ and $k_i>0$.
This prediction 
was shown in [GeP] to be implied by the
degree 0 Virasoro conjecture applied to $\proj^2$. 

\vspace{+10pt}
\noindent{\bf Acknowledgments.}
We thank
D. Zagier
for his proof of the Bernoulli identity 
required for (\ref{ihop})
which
resisted our best efforts. His argument appears in Section 4.4.
Conversations with E. Getzler, T. Graber, and E. Looijenga
played an important role in our work. 
This research was partly pursued at the Scuola Normale Superiore di
Pisa and Mathematisches Forschungsinstitut Oberwolfach in the summer of 1998.
The authors were partially supported by National Science
Foundation grants DMS-9801257 and DMS-9801574.

\section{\bf{Reconstruction equations}}

\subsection{Mumford's calculation}
We start by interpreting  Mumford's 
beautiful Grothendieck-Riemann-Roch calculation  [Mu]
in Gromov-Witten theory.
Let $\CM$ be a nonsingular variety (or Deligne-Mumford stack).
Let $\pi:\UM\rarr \CM$ be a flat family of genus $g$ pre-stable
curves (the fibers of $\pi$ are complete, connected, and
reduced, with
only nodal singularities). Assume the variation of this
family is maximal in the following sense:
the Kodaira-Spencer map
\begin{equation}
\label{asss} 
T\CM_m \rarr \text{Ext}^1(\Omega_{\UM_m}, \oh_{\UM_m})
\end{equation}
is surjective for every point $m\in M$.
In this case, the following facts are well-known
from the deformation theory of nodal curves:
\begin{enumerate}
\item[(i)] $\UM$ is a nonsingular variety (or Deligne-Mumford stack).
\item[(ii)] The singular locus of $\pi$ 
(the locus of nodes of the fibers) is a nonsingular
variety $Z$ of pure codimension $2$.
The image $\pi(Z)=\partial \CM$ is a divisor with normal crossings in $\CM$.
\item[(iii)] There is a natural \'etale double cover
$\epsilon:\tilde{Z} \rarr Z$ obtained from the 2-fold choice of
branches incident at the nodes corresponding to points of $Z$.
\item[(iv)] There are natural line bundles
$\LL$, $\overline{\LL}$ on $\tilde{Z}$
corresponding to the cotangent directions
along the branches.
\item[(v)] There is a canonical isomorphism  $\epsilon^*(\text{Nor}_{Z/\UM})
= \LL^* \oplus \overline{\LL}^*$. 
\end{enumerate}
Let $\iota: \tilde{Z} \rarr \CM$ denote the
natural composition. Let $\psi, \overline{\psi} \in A^1(\tilde{Z})$ denote
the first Chern classes of $\LL, \overline{\LL}$ respectively (Chow
groups will always be taken with $\Q$-coefficients).
The morphism $\iota$ is generically $2-1$ onto the
divisor $\partial \CM$. Let $\kappa_l=
\pi_*(c_1(\omega_\pi)^{l+1})\in
A^l(\CM)$.

The Hodge bundle
 is defined  on $\CM$ by $\hodge=\pi_* \omega_\pi$.
Mumford calculates $\ch (\hodge)$ in $A^*(\overline{M}_g)$ via the
Grothendieck-Riemann-Roch formula. 
As he uses only properties (i-v) above for the family
$\pi: \overline{M}_{g,1} \rarr \overline{M}_g$, 
his argument applies verbatim to the more general
setting considered here.

\vspace{+10pt}
\noindent{\bf Theorem} (Mumford).
$$\ch (\hodge) =
g + \sum_{l=1}^{\infty} \frac{ B_{2l}}{(2l)!}
\cdot \Big( \kappa_{2l-1} + \frac{1}{2} \iota_* \sum_{i=0}^{2l-2}
(-1)^i \psi^i \overline{\psi}^{2l-2-i} \Big)$$
in $A^*(\CM)$.

\vspace{+10pt}
\noindent The discrepancies between the above formula
and [Mu] are due to a differing Bernoulli number convention
and a typographic error in the $\kappa$ term of [Mu]. In our
formulas, the Bernoulli numbers are defined by:
\begin{equation*}
\dfrac{t}{e^t-1} = \sum\limits_{m=0}^{\infty}B_m\dfrac{t^m}{m!}.
\end{equation*}

\subsection{Gromov-Witten theory}
Let $X$ be a nonsingular projective variety over $\com$.
Let $\overline{M} =\overline{M}_{g,n}(X,\beta)$ be the moduli stack of
stable maps to $X$ representing the class $\beta \in H_2(X, \Z)$.
Let $[\overline{M}]^{vir} \in A_{*} (\overline{M})$
denote the virtual class in the expected dimension
[BF], [B], [LiT].
A direct analogue of Mumford's result holds for the
universal family over $\overline{M}$.

Virtual divisors in $\overline{M}$ are of two types.
First, stable splittings
\begin{equation}
\label{wqwt} 
\xi =(g_1+g_2=g, A_1 \scup A_2 =[n], \beta_1+\beta_2=\beta)
\end{equation}
index  virtual divisors in $\overline{M}$ corresponding
to maps with reducible domain curves. Define
\begin{equation}
\label{tyty}
\Delta_\xi = \overline{M}_{g_1, A_1+*}(X, \beta_1) \times_X
\overline{M}_{g_2, A_2+\bullet}(X, \beta_2) 
{\rarr} \overline{M}
\end{equation}
to
be the virtual divisor corresponding to the data $\xi$.
The fibered product in (\ref{tyty}) is taken with respect
to the evaluation maps $e_*, e_\bullet$
corresponding to the markings $*, \bullet$.
The virtual class of $\Delta_\xi$ is determined by:
$$[\Delta_\xi]^{vir} = [\overline{M}_{g_1, A_1+*}(X, \beta_1)]^{vir} 
\times [\overline{M}_{g_2, A_2+\bullet}(X, \beta_2)]^{vir}
\scap (e_*\times e_\bullet)^{-1} (\delta)$$
where $\delta \subset X\times X$ is the diagonal
(this is Axiom 4 of [BM]).

For $g\geq 1$, there is an additional
virtual divisor $\Delta_0$ corresponding to irreducible
nodal domain curves:
$$\Delta_0 = \overline{M}_{g-1, [n]+\{*, \bullet\}}(X, \beta) \scap
(e_*\times e_\bullet)^{-1} (\delta) {\rarr} 
\overline{M}$$
where $\delta \subset X\times X$ is the diagonal.
By Axiom 4, 
$$[\Delta_0]^{vir} = [\overline{M}_{g-1, [n]+\{*,\bullet\}}(X, \beta)]^{vir} 
\scap (e_*\times e_\bullet)^{-1} (\delta).$$

Let $\Delta$ be the set of all {\em ordered} splittings (\ref{wqwt})
indexing reducible divisors (with repetition) union $\{ 0\}$
for the irreducible divisor. 
There is a natural map
$$\iota: \bigcup_{\xi \in \Delta} \Delta_\xi \rarr \overline{M}$$
where the domain is the disjoint union.

Consider the morphism:
$$\overline{M} \rarr \mathfrak{M}_g$$
where the right side is the Artin stack of pre-stable 
genus $g$ curves.
For $0\leq j \leq g$, let
$$B_j = \mathfrak{M}_{j,*} \times \mathfrak{M}_{g-j,\bullet}.$$
Let $B_{irr}= \mathfrak{M}_{g-1,\{*,\bullet\}}$.
These Artin stacks admit natural maps
$\nu_0, \ldots, \nu_g$, $\nu_{irr}$ to $\mathfrak{M}_g$.
Let $\Delta^j \subset \Delta$ be the
subset with (ordered) genus splitting $g_1=j$, $g_2=g-j$.
Let $\Delta^{irr}= \{0\}$.
Certainly, 
$$\bigcup_{\xi\in \Delta^j} \Delta_\xi \eqq  B_j\times
_{\mathfrak{M}_g} \overline{M} $$
for $j\in \{0, \ldots,g, irr\}$.
The Isogeny Axiom of [BM] implies for each such $j$,
\begin{equation}
\label{ewee}
\nu_j^{!} [\overline{M}]^{vir} = \sum_{\xi\in \Delta^j}
[\Delta_\xi]^{vir}.
\end{equation}
This is one of the most important properties of the virtual class.

\begin{pr} \label{hodg}
$$ \ch(\hodge) \scap [\overline{M}]^{vir} = g [\overline{M}]^{vir}$$
$$+ \sum_{l=1}^{\infty} \frac{ B_{2l}}{(2l)!}
\cdot \Big( \kappa_{2l-1} \scap [\overline{M}]^{vir}
+ \frac{1}{2} \iota_*\sum_{\xi \in \Delta}
\sum_{i=0}^{2l-2}
(-1)^i \psi_*^i \psi_\bullet^{2l-2-i} \scap [\Delta_\xi]^{vir}\Big)$$
in $A_*(\overline{M})$.
\end{pr}

\bpf We will find a  
nonsingular
Deligne-Mumford stack $\CM$
with a family of curves
$\pi:\UM \rarr \CM$ satisfying assumption (\ref{asss}) and
an embedding 
$$\overline{M} \rarr \CM$$
such that $\UM$ restricts to the universal family over
$\overline{M}$:
\begin{equation*}
\begin{CD}
 U 
@>>>  \UM \\
@V{\pi}VV   @V{\pi}VV \\
 \overline{M} @>>> \CM.
\end{CD}
\end{equation*}
Following the notation of Section 1.1, we see
$$\tilde{Z}= \bigcup_{j\in\{irr,0,\ldots,g\}} B_j \times_{\mathfrak{M}_g}
\CM,$$
$$\tilde{Z} \times_\CM \overline{M} = 
\bigcup_{j\in\{irr,0,\ldots,g\}}
B_j \times_{\mathfrak{M}_g} \overline{M}.$$
We may then apply Mumford's Theorem to the map
$\pi: \UM\rarr \CM$. 
Intersecting Mumford's formula with
$[\overline{M}]^{vir}$ yields:
$$ \ch(\hodge) \scap [\overline{M}]^{vir} = g [\overline{M}]^{vir}$$
$$+ \sum_{l=1}^{\infty} \frac{ B_{2l}}{(2l)!}
\cdot \Big( \kappa_{2l-1} \scap [\overline{M}]^{vir}
 +  
{\frac{1}{2}}\iota_* \sum_{j\in\{irr,0,\ldots,g\}}
 \sum_{i=0}^{2l-2}
(-1)^i \psi_*^i \psi_\bullet^{2l-2-i}
\scap \nu_j^{!} [\overline{M}]^{vir}
\Big)$$
in $A_*(\overline{M})$.
The proposition then follows immediately from (\ref{ewee}).

The construction of the required family $\pi: \UM \rarr \CM$
starts with a general observation.
Let 
\begin{equation}
\label{wef}
S\subset \proj^r \times B \rarr B
\end{equation} 
be a projective flat family of  genus $g$, degree $d$ pre-stable
curves over a quasi-projective base scheme $B$.
We show how to embed (\ref{wef}) in a family of curves
over a nonsingular base satisfying assumption (\ref{asss}).

Let $\mathcal{L}=\oh_{\proj^r}(1)$.
By standard boundedness arguments, there exists
an integer $\alpha$ satisfying: 
\begin{equation}
\label{sdsdss}
H^1(S_b, \mathcal{L}^\alpha_b)=0
\end{equation}
for all $b\in B$. Consider the Veronese embedding
$$\proj^n \rarr \proj^{\binom{n+\alpha}{\alpha}-1}$$
Then, there is a canonical map
$$\phi_1: B \rarr \mathcal{H},$$
where
$\mathcal{H}$ is the Hilbert scheme 
of genus $g$, degree $\alpha d$ curves in 
$\proj^{\binom{n+\alpha}{\alpha}-1}$.
The vanishing (\ref{sdsdss}) easily implies $\mathcal{H}$
is nonsingular of expected dimension in a Zariski open set
$\mathcal{H}^0 \subset \mathcal{H}$
containing 
$\text{Im}(\phi_1)$.
Assumption (\ref{asss}) for the
universal family $\mathcal{U}^0 \rarr \mathcal{H}^0$
also is a direct consequence of (\ref{sdsdss}).
Let $\phi_2:B \rarr X$ be a closed embedding in a 
nonsingular scheme $X$.
Finally, the diagram
\begin{equation*}
\begin{CD}
S
@>>>  \mathcal{U}^0 \times X \\
@VVV   @VVV \\
B @>{\phi_1 \times \phi_2}>> \mathcal{H}^0 \times X.
\end{CD}
\end{equation*}
is the required construction 
for the given family
$S \rarr B$.

In [FP], $\overline{M}$ is constructed as a Deligne-Mumford
quotient stack $\text{Hilb}/\G$ of a reductive group action
on a Hilbert scheme of pointed graphs. 
The universal family $U \rarr \overline{M}$ is simply
the stack quotient of the universal family $\mathcal{U} \rarr
\text{Hilb}$.
The above construction applied $\G$-equivariantly
to $\mathcal{U} \rarr \text{Hilb}$ directly yields
the required construction for the Proposition
(see also [GrP] where embeddings of $\overline{M}$ in
nonsingular Deligne-Mumford stacks are constructed).
\epf

\subsection{Theorem 1}

Let $X$ be a nonsingular projective variety of dimension $r$.  
Let $\gamma_0, \ldots, \gamma_N$ be a graded $\Q$-basis of
$H^*(X, \Q)$. We take $\gamma_0$ to be the identity element.
Let $g_{ef}= \int_X \gamma_e\scup \gamma_f$, and
let $g^{ef}$ be the inverse matrix.
The 
descendent
Gromov-Witten invariants of $X$
$$\lan  \prod_{i=1}^n \tau_{k_i}(\gamma_{a_i}) \ran_{g,\beta}^X =
\int_{[\mgn(X,\beta)]^{vir}} 
\prod_{i=1}^n \psi_i^{k_i} \scup e_i^*(\gamma_{a_i})$$
may be organized in a generating function
\begin{equation*}
F^X = \sum_{g\ge 0} \hbar^{g-1} F^X_g,
\end{equation*}
where
$$
F_g^X = \sum_{\beta \in H_2(X,\Z)} q^\beta \sum_{n\ge 0} \frac{1}{n!}
\sum_{\substack{a_1\dots a_n \\ k_1 \dots k_n}} t_{k_n}^{a_n} \dots
t_{k_1}^{a_1} \lan \tau_{k_1}(\gamma_{a_1}) \dots \tau_{k_n}(\gamma_{a_n})
\ran_{g,\beta}^X. 
$$

We introduce here an analogous
generating function $F^X_{\hodge}$ for the Hodge integrals
over moduli stacks of maps to $X$.
For each odd positive integer, let the
variable $s_i$ correspond to $\ch_i (\hodge)$. By Mumford's
relations [Mu], the even components of $\ch(\hodge)$
vanish (for all genera).
Let $t,s$ denote the sets of variables $\{t_i^j\}$, $\{s_i\}$
respectively.
The Hodge integrals
$$\lan  \prod_{i=1}^n \tau_{k_i}(\gamma_{\alpha_i}) 
\prod_{j=1}^{m} \ch_{b_j}(\hodge)  \ran_{g,\beta}^X =$$
$$\int_{[\mgn(X,\beta)]^{vir}} 
\prod_{i=1}^n \psi_i^{k_i} \scup e_i^*(\gamma_{\alpha_i})
\ \scup \prod_{j=1}^{m}  \ch_{b_j}(\hodge)$$
define formal functions
$$F^X_{g,\hodge} (t,s) =$$ 
$$\sum_{\beta \in H_2(X,\Z)} q^\beta \sum_{n,m\ge 0} \frac{1}{n!m!}
\sum_{\substack{k_1 \dots k_n \\ a_1 \dots a_n
\\ b_1  \dots b_m}} 
\prod_{i=1}^n
t_{k_i}^{a_i}
\prod_{j=1}^m s_{b_j} \lan 
\prod_{i=1}^n\tau_{k_i}(\gamma_{a_i}) \prod_{j=1}^m \ch_{b_j}(\hodge)
\ran_{g,\beta}^X.$$
As before, we define
$F^X_\hodge = \sum_{g\ge 0} \hbar^{g-1} F^X_{g, \hodge}.$
This function is related to the descendent generating function
by restriction: $F^X_\hodge | _{s=0} = F^X$.
Finally, 
let $Z^X_\hodge= \text{exp}(F^X_\hodge)$.

Formulas involving the cotangent line classes 
and the Chern character of the Hodge bundle
yield
the following consequence of Proposition \ref{hodg}.
For $l \geq 1$, define a formal differential
operator:
$$D_{2l-1} = $$
$$-\frac{\partial}{\partial{s_{2l-1}}}  
+\frac{B_{2l}}{(2l)!} \Big(\frac{\partial}{\partial {t^0_{2l}}} 
- \sum_{i=0}^{\infty} \sum_{j=0}^{N}
t_i^j \frac{\partial}{\partial t_{i+2l-1}^{j}}
+\frac{\hbar}{2}
\sum_{i=0}^{2l-2} (-1)^i g^{ef} \frac{\partial}{\partial t_{i}^e }
\frac{\partial}{\partial t_{2l-2-i}^f}\Big),$$
as usual the sum over the indices $e$,$f$ is suppressed.
\begin{pr} \label{diff}
For all $l\geq 1$,
$D_{2l-1} Z^X_\hodge =0.$
\end{pr}
\bpf
Let $\overline{M}= \overline{M}_{g,n}(X,\beta)$ as in
Section 1.1. 
Three formulas are needed to deduce this vanishing from
Proposition \ref{hodg}.

Let $d$ be the virtual dimension of $\overline{M}$.
The Chow class $\kappa_{2l-1} \scap [\overline{M}]^{vir}$
has dimension $d-2l+1$.
The first formula is:
\begin{equation}
\label{dffd}
\prod_{i=1}^n \psi_i^{k_i}\scup e_i^*(\gamma_{a_i})
   \ \scap (\kappa_{2l-1} \scap [\overline{M}]^{vir})=
\end{equation}
$$
\lan \tau_{2l}(\gamma_0) \prod_{i=1}^n \tau_{k_i}(\gamma_{a_i}) 
\ran^X_{g,\beta} -
\sum_{i=1}^n
\lan \tau_{k_i+2l-1}(\gamma_{a_i}) \prod_{j\neq i} \tau_{k_j}(\gamma_{a_j}) 
\ran^X_{g,\beta},$$
where the cohomology product on the left side has codimension
$d-2l+1$.
It 
follows from viewing
the universal family over $\overline{M}$ as
$\overline{M}_{g,n+1}(X, \beta)$ and applying
the standard comparison
results for cotangent lines (see [W]). The only virtual class
property needed is the equality 
$$[\overline{M}_{g,n+1}(X,\beta)]^{vir} = \pi^* [\overline{M}]^{vir}$$
which is an Axiom in [BM].

The second and third required formulas address the
behavior of the Chern character of the Hodge bundle
when restricted to the virtual boundary divisors.
Let $\xi\in \Delta$ correspond to a virtual boundary
divisor with genus splitting $g_1+g_2=g$.
Let $\hodge_g$ denote the Hodge bundle on $\overline{M}$.
Let $\hodge_{g_1}$, $\hodge_{g_2}$ denote the Hodge
bundles obtained from the two factors in (\ref{tyty}).
The natural restriction sequence on $\Delta_\xi$:
$$0 \rarr \hodge_{g_1} \rarr \iota^* \hodge_{g} \rarr
\hodge_{g_2} \rarr 0$$
implies the formula
\begin{equation}
\label{ssec}
\ch(\hodge_{g_1}) + \ch(\hodge_{g_2}) = \iota^* \ch(\hodge_g)
\in A^*(\Delta_\xi).
\end{equation}
Similarly, for the irreducible virtual divisor $\Delta_{0}$,
the residue sequence
$$0 \rarr \hodge_{g-1} \rarr \iota^* \hodge_g \rarr \oh_{\Delta_0}
\rarr 0$$
implies the formula
\begin{equation}
\label{ttec}
\ch(\hodge_{g-1}) = \iota^* \ch(\hodge_g)
\in A^*(\Delta_0).
\end{equation}
Proposition 2 is a formal consequence of Proposition 1
and equations (\ref{dffd}-\ref{ttec}).
\epf

The generating function $F^X_{\hodge}$ is determined
by the initial $s=0$ conditions (specified by $F^X$) and
the differential equations from Proposition \ref{diff}. Thus,
Theorem 1 is proven. 

We end this section with some remarks following from
Proposition \ref{diff}. All the Chern classes
of the Hodge bundle vanish in genus 0. 
Hence, $\partial F^X_{0,\hodge}/ \partial 
s_{2l-1}=0$. The vanishing $D_{2l-1} Z^X_\hodge =0$ 
analyzed at order $\hbar^{-1}$ then yields 
universal relations among genus 0 descendent invariants
of $X$. The relation obtained for $l=1$ is coincides precisely with
a derivative of
$\tilde{L}_1$ (defined in [EHX] and used in
the proof of the genus 0 Virasoro
constraints). Proposition \ref{diff} also
yields geometric interpretations of several related
equations in the latter proof (see [Ge]).

In fact, Proposition \ref{diff} yields many more new universal relations
among pure descendent invariants.
For example, the classes $\ch_{2l-1}(\hodge)$ vanish
in $A^*(\overline{M}_g)$ for $l>g$.
Hence, generalizations
of the above $g=0$ equations to higher genus
are obtained from 
$$\frac{\partial F_{g,\hodge}^X}{\partial s_{2l-1}}=0 \ \ \ (l>g),$$
and the vanishing at order $\hbar^{g-1}$ in
$D_{2l-1} Z^X_\hodge =0$. 
The resulting relation is an efficient topological
recursion relation (TRR) for $\tau_{2l}$ in genus $g<l$.
Note the Bernoulli number drops out of these relations.

A more sophisticated method of obtaining pure
descendent equations from Proposition \ref{diff} is to
construct combinations of the operators $D_{2l-1}$
that serve to introduce the Chern classes of
$\hodge$. The Chern classes of $\hodge$ certainly vanish in
degrees greater than $g$ on $\overline{M}_g$.
One obtains from Proposition \ref{diff} relations in
degree greater than $g$ (for each $g$). It
would be interesting to understand these
equations and their relation to TRR and 
the Virasoro constraints even in the point case.

Finally, while the Hodge integrals
\begin{equation*}
\int_{[\mgn(X,\beta)]^{vir}} 
\prod_{i=1}^n \psi_i^{a_i} \scup e_i^*(\gamma_i)
\scup \prod_{j=1}^g \lambda_j^{b_j}
\end{equation*} 
are determined by Proposition \ref{diff} and $F^X$, the
relations satisfied by the natural generating functions
of these integrals do not appear easy to write.

\section{\bf{Relations via virtual localization}} 
\subsection{The localization formula}
We review here the virtual localization formula
in Gromov-Witten theory [GrP] in the special case of
degree 1 maps to $\proj^1$. While our strategy for 
obtaining relations among Hodge integrals may be pursued
in much greater generality, only this special case is
required for Theorem 2.

Let $\proj^1=\proj(V)$ where $V=\com \oplus \com$.
Let $\com^*$ act diagonally on $V$:
\begin{equation}
\label{repp}
\xi\cdot (v_1,v_2) = ( v_1, 
\xi \cdot v_2).
\end{equation}
Let $p_1, p_2$ be the fixed points of the corresponding
action on $\proj(V)$.
An equivariant lifting  of $\com^*$ to a line bundle $L$
over 
$\proj(V)$ is uniquely determined by the weights $[l_1,l_2]$
of the fiber
representations at the fixed points 
$$L_1= L|_{p_1}, \ \ \ L_2= L|_{p_2}.$$
The canonical lifting of $\com^*$ to the
tangent bundle, $\text{Tan}$, has weights $[1,-1]$.
There is a scaling lifting of $\com^*$ to $\oh_{\proj(V)}$ for
each integer $\alpha$ with weights $[\alpha, \alpha]$.
For each integer $\beta$, there is a $\com^*$-lifting
to $\oh_{\proj(V)}(-1)$ with weights $[\beta, \beta+1]$.

Let $g\geq 1$.
Let $\mgp= \overline{M}_{g,0}(\proj(V),1)$ be the
moduli stack of stable, genus $g$, unpointed maps
to $\proj(V)$ of degree 1. Let 
\begin{equation}
\label{zzzz}
\pi: U \rarr \mgp , \ \ \ \mu: U \rarr \proj(V)
\end{equation}
be the universal curve and universal map over the
moduli stack. The representation (\ref{repp}) canonically
induces $\com^*$-actions on $U$ and $\mgp$ compatible
with the maps $\pi$ and $\mu$.

The virtual dimension of $\mgp$ is $2g$. 
There are two natural rank $g$ bundles on $\mgp$:
$R^1\pi_*(\mu^*{\oh_{\proj(V)}})$ and 
$R^1\pi_* (\mu^* \oh_{\proj(V)}(-1))$.
Let $x, y$ denote the respective top Chern
classes of these bundles in $A^g(\mgp)$.
The
following two integrals against the virtual
class $[\mgp]^{vir} \in A_{2g}(\mgp)$ will be considered:
\begin{equation}
\label{rrrr}
\int_{[\mgp]^{vir}} x \scup y,
\ \
\int_{[\mgp]^{vir}} y \scup y.
\end{equation}
The virtual localization formula will be used to
compute these integrals with respect to various
linearizations on $\oh_{\proj(V)}$ and $\oh_{\proj(V)}(-1)$.

The fixed locus $X$ of the $\com^*$-action on
$\mgp$ is a disjoint union of irreducible
components 
$$X=\bigcup_{\substack{g_1+g_2=g\\g_1,\,g_2 \geq 0}} X_{g_1,g_2}.$$
The component $X_{g_1,g_2}$ corresponds to
the loci of maps where subcurves of genus $g_1$ and
$g_2$ are contracted to the fixed points $p_1$ and
$p_2$ respectively. The fixed locus is naturally isomorphic
to $\overline{M}_{g_1,1} \times \overline{M}_{g_2,1}$ (where
$\overline{M}_{0,1}$ is defined to be a point).
Moreover, the induced fixed stack structure   
on $X_{g_1, g_2}$ is simply the reduced nonsingular structure [GrP].
The cotangent line and $\lambda$ classes of the
two factors yield cohomology classes on $X_{g_1, g_2}$
via pull-back. Let $\psi_1, \psi_2$ denote the 
cotangent line classes from the factors
$\overline{M}_{g_1,1}$ and $\overline{M}_{g_2,1}$ respectively.
For $k\in \Z$, let
$$\HH_1(k)= \sum_{i=0}^{g_1} k^i \lambda_{g_1-i} \in 
A^*(\overline{M}_{g_1,1}),$$
$$\HH_2(k)= \sum_{i=0}^{g_2} k^i \lambda_{g_2-i} \in 
A^*(\overline{M}_{g_2,1}).$$
We note Mumford's formula $c(\hodge)\cdot c(\hodge^*)=1$ 
implies
\begin{equation}
\label{mmumm}
\HH_i(-1) \HH_i(1) = (-1)^{g_i},
\end{equation}
$$ \HH_i(0) \HH_i(0) = \delta_{g_i 0}.$$
These sums $\HH_i(k)$ will be convenient for the formulas below.

Let $\iota: X \rarr \mgp$ be the inclusion.
The virtual localization formula is:
\begin{equation}
\label{wwww}
\iota_* \sum_{g_1+g_2=g} \frac{[X_{g_1,g_2}]}
{c_{\text{top}}(\text{Nor}^{vir}_{g_1,g_2})} = [\mgp]^{vir} \ \in 
H^*_{\com^*} (\mgp)[1/t].
\end{equation}
The virtual normal bundle $\text{Nor}^{vir}_{g_1,g_2}$ is 
isomorphic 
in equivariant $K$-theory on $X_{g_1,g_2}$ to the sum:
$$[\psi_1 \otimes \text{Tan}_1] + [\psi_2
\otimes \text{Tan}_2] +
[\pi_*\text{Tan}]-[R^1\pi_*\text{Tan}] -[\text{Aut}]$$
(see [GrP]).
Let $\gamma\in H^{4g}_{\com^*}(\mgp)$. 
After an expansion of the virtual normal
contribution, equation (\ref{wwww})  yields 
an explicit integration formula for $\gamma$:
\begin{equation}
\label{xplict}
\int_{[\mgp]^{vir}} \gamma = \sum_{g_1+g_2=g} \int_{X_{g_1,g_2}}
(-1)^{g} \ \iota^*(\gamma) \ \frac {\HH_1(-1)}
{1-\psi_1}   \ \frac{\HH_2(1)}{1+\psi_2}.
\end{equation}

\subsection{Relations}
Application of formula (\ref{xplict}) to the integrals
(\ref{rrrr}) yields the following linearization dependent
equations.
We find
$$
\int_{[\mgp]^{vir}} x \scup y = (-1)^g I_g(\alpha, \beta)$$
with respect to the linearizations 
$[\alpha, \alpha]$ on $\oh_{\proj(V)}$ and $[\beta, \beta+1]$ on 
$\oh_{\proj(V)}(-1)$
where
\begin{equation}
\label{iiii}
I_g(\alpha, \beta) =
\sum  \int_{X_{g_1,g_2}}
 \frac {\HH_1(-1) \HH_1(-\alpha) \HH_1(-\beta)}
{1-\psi_1}   \frac{\HH_2(-1)\HH_2(\alpha) \HH_2(\beta+1)}{1-\psi_2}.
\end{equation}
Similarly, 
$$
\int_{[\mgp]^{vir}} y \scup y = (-1)^g J_g(\alpha, \beta)$$
with respect to the linearizations $[\alpha, \alpha+1]$,
$[\beta, \beta+1]$ on the two copies of $\oh_{\proj(V)}(-1)$ where
\begin{equation}
\label{jjjj}
J_g(\alpha, \beta) =
\sum \int_{X_{g_1,g_2}}
\frac {\HH_1(-1) \HH_1(-\alpha) \HH_1(-\beta)}
{1-\psi_1}   \frac{\HH_2(-1)\HH_2(\alpha+1)\HH_2(\beta+1)}{1-\psi_2}.
\end{equation}
Hence, we have obtained the relations
\begin{equation}
\label{rell}
I_g(\alpha, \beta)= I_g(\alpha', \beta'), \ \
J_g(\alpha, \beta)= J_g(\alpha', \beta')
\end{equation}
for all integers $\alpha, \alpha', \beta, \beta'$.

For $\xi\in \Z$, define the series $f_\xi(t)\in \Q[[t]]$ by:
$$f_\xi(t)=  1+ \sum_{g\geq 1} t^{2g} \int_{\overline{M}_{g,1}}
\frac{\HH(\xi)}{1-\psi_1} =
1 + \sum_{g\geq 1} 
\sum_{i=0}^{g} t^{2g} \xi^i
\int_{\overline{M}_{g,1}} \psi_1^{2g-2+i} \lambda_{g-i}.$$

\begin{pr} \label{fxi}
For $\xi\in \Z$, $f_\xi(t)= f_0(t)^{\xi+1}$.
\end{pr}
\bpf 
By the integration formulas (\ref{iiii}-\ref{jjjj})
together with Mumford's relations
(\ref{mmumm}), we find:
$$ 1+ \sum_{g\geq 1} t^{2g} I_g(0,0) = f_0(it),$$
$$ 1+ \sum_{g \geq 1} t^{2g} J_g(0,-1)= f^2_0(it).$$
We will consider the relations:
$$ 1+ \sum_{g\geq 1} t^{2g} I_g(\xi,0) = f_0(it),$$
$$ 1+ \sum_{g \geq 1} t^{2g} J_g(0,\xi) = f^2_0(it).$$
Define a new series for $\xi\in \Z$:
$$ g_\xi(t)= 1 + \sum_{g\geq 1} t^{2g} \int_{\overline{M}_{g,1}}
\frac{\HH(-1) \HH(0) \HH(-\xi)}{1-\psi_1}.$$
The integration formulas imply:
$$   1+ \sum_{g\geq 1} t^{2g} I_g(\xi,0) = g_\xi(t) f_\xi(it),$$
$$1+ \sum_{g\geq 1} t^{2g} J_g(0,\xi) =   g_\xi(t) f_{\xi+1}(it).$$
We then deduce the equations:
$$ g_\xi(t) f_\xi(it)= f_0(it), \ \ 
g_\xi(t) f_{\xi+1}(it)= f_0^2(it).$$
Hence, $f_{\xi+1}(it)= f_\xi(it) f_0(it)$ for
all $\xi \in \Z$.
The proposition now follows easily 
by induction (as it is true for $\xi=0$).
\epf

In order to determine the functions $f_\xi(t)$, it suffices
to compute only $f_{-2}(t)= f_0(t)^{-1}$. This
calculation too may be accomplished via localization
relations, but a shorter and more elegant derivation
by classical curve theory will be given in Proposition \ref{hyp}.

To show the flavor of Hodge relations obtained from
localization, we mention two further examples.
The formula:
\begin{equation}
\label{wk}
1 + \sum_{g\geq 1} t^g \int_{\overline{M}_{g,1}} \psi^{3g-2} = 
\text{exp}(t/24)
\end{equation}
is a well known consequence of Witten's conjecture (Kontsevich's
theorem). It is a nice exercise to prove this formula
via Hodge relations obtained from localization on the
stack of maps to $\proj^1$. A geometric proof of
(\ref{wk}) will be given in the next section.

Let $\gamma \in H^2(\proj^1)$ be the point class.
The integral
\begin{equation}
\label{asas}
\int_{[\overline{M}_{g,1}(\proj^1,d)]^{vir}}
x\scup y \scup e_1^*(\gamma^{d})
\end{equation}
clearly vanishes for $d\geq 2$ 
(as before $x$ and $y$ are the top Chern classes
of the vector bundles obtained from the
higher direct images of $\mu^*(\oh_{\proj(V)})$
and $\mu^*(\oh_{\proj(V)}(-1))$ respectively, $e_1$
is the evaluation map corresponding to the marking).
When (\ref{asas}) is computed by localization
with an appropriate choice of linearization, the
following Hodge relation is found:
$$ \sum_{ m\in \text{Part}(d)}
\frac{(-1)^{d+l(m)}  \prod_i m_i^{m_i}}
{\text{Aut}(m) \prod_i m_i  \prod_i m_i!} 
\ \int_{\overline{M}_{g,l(m)+1}} \frac{\lambda_g}{\prod_i(1-m_i \psi_i)}=0$$
where $m=\{m_1, \ldots, m_{l(m)}\}$ is a partition
of $d$.
We have checked algebraically that the Virasoro prediction
(\ref{lamg}) of [GeP] satisfies these relations. As yet, we are unable
to prove (\ref{lamg}) via Hodge relations of this type.

\section{\bf{Relations via classical curve theory}}
\subsection{Relations via the canonical system}
In this section, we derive several relations among Hodge integrals
from classical curve theory. The starting point is [Mu].
The base-point-freeness of the canonical system on a smooth curve
can be formulated as the surjectivity of the natural map
$\pi^*\hodge\to\LL_1$ on $C_g=M_{g,1}$. This gives rise to an exact
sequence
$$0\to F\to\pi^*\hodge\to\LL_1\to0$$
with $F$ locally free of rank $g-1$.
Hence one finds on $M_{g,1}$ the relations
$$\left(\frac{c(\hodge)}{1+\psi_1}\right)_j=0\qquad(j\ge g).$$
If we want to extend these relations to 
$\overline{M}_{g,1}$,
we must take into account the stable pointed curves
for which $\LL_1$ is not generated by global sections. As Mumford observes,
the global sections generate the subsheaf of $\LL_1$ that is zero at
all disconnecting nodes and on all smooth rational curves all of whose nodes
are disconnecting. Let us denote for $2\le i\le g$ by $X_i$ the locus
of stable one-pointed curves of genus $g$ consisting of a stable
$(i+1)$-pointed rational curve with $i$ tails (stable one-pointed curves
of positive genus; the $i$ genera sum to $g$) attached to the last
$i$ marked points. It follows that the relations above hold on
$\overline{M}_{g,1}$ modulo a class supported on the loci 
$X_2, \dots, X_g$. (Note that $X_2$ is the locus of disconnecting nodes
in the universal curve.)

Since the moduli stack of $(i+1)$-pointed rational curves has dimension
$i-2$ we have that $\psi_1^{i-1}$ is 0 on $X_i$. Hence $\psi_1^{g-1}$
is 0 on all these loci; we find the relations
$$\left(\frac{c(\hodge)}{1+\psi_1}\right)_j=0\qquad(j\ge 2g-1)$$
on $\overline{M}_{g,1}$. For $j=3g-2$,
we find 
\begin{equation}
\label{zeroz}
 \int_{\overline{M}_{g,1}} \frac{\HH(1)}{1+\psi_1} = 0
\end{equation}
(in the notation of Section 2).
This identity implies $f_{-1}(t) = 1$  which is also a
consequence of Proposition \ref{fxi}.

If instead we intersect the relation for $j=g$ with $\psi_1^{g-2}$, we
find
\begin{equation} \label{xg*}
\left(\frac{c(\hodge)}{1+\psi_1}\right)_{2g-2}=*\psi_1^{g-2}[X_g]_Q\,.
\end{equation}
Here $[\hphantom{x}]_Q$ denotes the $Q$-class or fundamental class in the
sense of stacks as in [Mu]. The coefficient $*$ can be determined
by intersecting with the locus $Y$ parametrizing one-pointed
irreducible curves with $g$ nodes (hence with rational normalization)
and their degenerations. Let $Z=X_g\cap Y$; this is the locus
of one-pointed curves consisting of a stable $(g+1)$-pointed
rational curve with $g$ singular elliptic tails attached. The intersection
is transverse in the universal deformation space, so that
$[X_g]_Q\cdot[Y]_Q=[Z]_Q$; it is easy to see that $\psi_1^{g-2}$
times this class equals $\dfrac1{2^gg!}\,$.

As the restriction of $\hodge$ to $Y$ is trivial, the
intersection of the left side of
(\ref{xg*}) with $[Y]_Q$ is $\psi_1^{2g-2}[Y]_Q$.
This product  evaluates
to $\dfrac1{2^gg!}$ as well, since the natural map $\overline{M}_{0,2g+1}
\to Y$ has degree $2^gg!$. We conclude that the coefficient $*$
in (\ref{xg*}) is equal to $1$:
\begin{equation} \label{xg}
\left(\frac{c(\hodge)}{1+\psi_1}\right)_{2g-2}=
\psi_1^{g-2}(\psi_1^g-\psi_1^{g-1}\lambda_1+\dots+(-1)^g\lambda_g)=
\psi_1^{g-2}[X_g]_Q\,.
\end{equation}
Intersecting this relation with 
$\psi_1^g+\psi_1^{g-1}\lambda_1+\dots+\lambda_g$ gives just
$\psi_1^{3g-2}$ on the left  side, since $c(\hodge)
\cdot c(\hodge^*)=1$.
On the right side we obtain $\lambda_g\psi_1^{g-2}[X_g]_Q$ which
easily  evaluates
to $1/(24^gg!)$. We find another proof of the identity (\ref{wk}),
$$\int_{\overline{M}_{g,1}} \psi_1^{3g-2}=\frac1{24^gg!}.$$

\subsection{Relations via Weierstrass loci}
Above, our starting point was the base-point-freeness of the canonical system
on a smooth curve. We then extended some of the relations so obtained to the
moduli stack of stable curves. Below, we study hyperelliptic Weierstrass
points; this may be viewed as a first step in analyzing the very-ampleness
of the canonical system. We obtain the following result.

\begin{pr} \label{hyp}
$$
f_{-2}(t)=1+\sum_{g\geq 1} t^{2g}
\int_{\overline{M}_{g,1}}\psi_1^{2g-2}(\lambda_g-2\psi_1\lambda_{g-1}
+\dots+(-2\psi_1)^g)=\frac{\sin(t/2)}{t/2}\,.
$$
\end{pr}

\bpf
In [Mu] Mumford computed the class in $C_g$ of the locus
$WH_g$ of hyperelliptic Weierstrass points:
\begin{eqnarray*}
&&[WH_g]_Q = \left(c(\hodge^*)
\frac1{1-\psi_1}\frac1{1-2\psi_1}\right)_{g-1} \\
&&\qquad
= (2^g-1)\psi_1^{g-1}-(2^{g-1}-1)\psi_1^{g-2}\lambda_1+\dots
+(-1)^{g-1}(2^1-1)\lambda_{g-1}\,.\\
\end{eqnarray*}
Hence,
\begin{eqnarray*}
\psi_1[WH_g]_Q &=& \left( (2\psi_1)^g-\lambda_1(2\psi_1)^{g-1}+\dots
+(-1)^g\lambda_g\right) \\
&&-\left( (\psi_1)^g-\lambda_1(\psi_1)^{g-1}+\dots 
+(-1)^g\lambda_g\right). \\
\end{eqnarray*}
Let us suppose this identity continues to hold on ${\overline{C}}_g
=\overline{M}_{g,1}$
modulo classes on which $\psi_1^{2g-2}$ is zero. Then,
by the vanishing (\ref{zeroz}), the formula
for $f_{-2}(t)$ is equivalent to
\begin{equation}
\label{nnn}
\psi_1^{2g-1}[{\overline{WH}}_g]=\frac1{2^{2g-1}(2g+1)!}\,.
\end{equation}
Note 
the usual fundamental class appears on the left in (\ref{nnn}), 
as this is more
convenient in the sequel.
We will first prove identity (\ref{nnn}) and then verify the 
required assumption.

The space ${\overline{M}}_{0,2g+2}$ may be viewed as the moduli space
of stable hyperelliptic curves of genus $g$ with an ordering of the
$2g+2$ Weierstrass points. (The hyperelliptic automorphism is lost
in this identification, however.) The universal (ordered) hyperelliptic
curve is a double cover of ${\overline{M}}_{0,2g+3}$ (the universal curve
over ${\overline{M}}_{0,2g+2}$). The ramification locus is
${\overline{WH}}_g$ (ordered); the branch locus $B$ is 
$$\sum_{j=1}^{2g+2}\,D_{j,2g+3}\,,$$
where $D_{j,2g+3}$ is the boundary divisor corresponding to the
partition $\{j,2g+3\}\cup\{j,2g+3\}^c$ (note that the $2g+2$ divisors
are disjoint).
The reason we can compute $\psi_1^{2g-1}[{\overline{WH}}_g]$
is that $\psi_1$ on the double cover is a pullback from 
${\overline{M}}_{0,2g+3}$. Denote the double cover map by $f$, then
$\psi_1=f^*(\psi_{2g+3}-B/2)$. This follows from the Riemann-Hurwitz
formula; note that $\psi_{2g+3}$ has degree $-2+(2g+2)=2g$ on the fibers
of the map to ${\overline{M}}_{0,2g+2}$. Hence
$$
\psi_1^{2g-1}[{\overline{WH}}_g]=f_*(\psi_1^{2g-1}[{\overline{WH}}_g])
=(\psi_{2g+3}-\frac12B)^{2g-1}B=(-\frac12)^{2g-1}B^{2g}\,.
$$
The last equality holds because $\psi_{2g+3}$ is zero on every component
of $B$. Now $B$ consists of $2g+2$ disjoint components, each
isomorphic to ${\overline{M}}_{0,2g+2}$; the restriction of $B$ to itself
is then $-\psi_*$ if $*$ is the marked point corresponding to the node.
Hence
$$
\psi_1^{2g-1}[{\overline{WH}}_g]=(2g+2)(\frac12)^{2g-1}\psi_*^{2g-1}
=\frac{2g+2}{2^{2g-1}}\,.
$$
This is the answer in the ordered case; the formula for the unordered
case follows immediately.

It remains to verify the assumption made: that Mumford's formula
for $[WH_g]_Q$ valid on $C_g$ holds on ${\overline{C}}_g$ after
multiplying by $\psi_1^{2g-1}$.
One may prove Mumford's formula by observing that the locus of
hyperelliptic Weierstrass points is the degeneracy locus
$\{\rk\,\phi_2\le1\}$,
where $\phi_2:\hodge\to{\mathbb{F}}_2$ is the natural evaluation map
from the Hodge bundle to the jet bundle
$\mathbb{F}_2$ whose fiber
at $[C,p]$ is the vector space $H^0(C,K/K(-2p))$ of dimension 2.
The class of the locus is then given by Porteous's formula. 

In order to verify the assumption,
we must analyze the 
irreducible loci $\{L_i\}$ of singular stable pointed curves
included in the degeneracy locus 
$\{\rk\,\phi_2\le1\}$ and
show $\psi_1^{2g-1} [L_i]=0$.
If $[C,p]$ lies in the degeneracy locus,
it is easy to see
one of the following two possibilities must be satisfied:
\begin{enumerate}
\item[(a)]  $p$ lies on a nonsingular rational component $X$;
\item[(b)] $p$ is a
hyperelliptic Weierstrass point: the component $X$ containing $p$
is a possibly nodal hyperelliptic curve (of arithmetic
genus $h\geq 1$), and the point $p$
is a Weierstrass point on $X$. 
\end{enumerate}
Let $L_i$ be an irreducible boundary component of the
degeneracy locus. If (a) holds generically on $L_i$,
naive estimates show 
the moduli of the component $X$ (with marked nodes
and point $p$) is bounded by $(3g-4)/2$ parameters. 
Hence, $\psi_1^{2g-1}[L_i]$
is certainly 0 is this case. 
Suppose (b) holds generically
on $L_i$. 
We may assume  the generic total curve $C$
is not hyperelliptic, otherwise $L_i$ lies in the closure of $WH_g$
and is of dimension less than $2g-1$.
In particular, $C$ must be reducible. We will show  the marked
component $X$ has fewer than $2g-1$ moduli.
We may assume $X$ is nonsingular and
meets the rest of the curve in $m$ points. We have to show that
$2h-1+m<2g-1$. Since $h<g$ this is clear when $m=1$. When $m=2$,
$h=g-1$ doesn't result in a stable curve of genus $g$, so we are done.
For $m\ge3$, the maximal $h$ is obtained when rational curves are
attached. But attaching a $k$-pointed rational curve lowers $2h-1+m$
by $k-2$, so $2h-1+m$ is always smaller than $2g-1$.
This finishes the proof of Proposition \ref{hyp}.\epf

\subsection{Proof of Theorem 2}
Define the series $F(t,k)\in \Q[k][[t]]$ by
$$F(t,k)= 1+ \sum_{g\geq 1} \sum_{i=0}^{g} t^{2g} k^i
\int_{\overline{M}_{g,1}} \psi_1^{2g-2+i} \lambda_{g-i}.$$
By Propositions \ref{fxi} and \ref{hyp}
$$F(t,\xi)= f_\xi(t)=\Big( \frac{t/2}{\sin(t/2)} \Big)^{\xi+1}$$
for all $\xi \in \Z$. The equality of formal series
$$F(t,k)= \Big( \frac{t/2}{\sin(t/2)} \Big) ^{k+1}$$
then follows  immediately. Theorem 2 is proven. \qed

\section{\bf{Bernoulli identities and Theorems 3-4}}

\subsection{Proof of Theorem 3}
Let $\overline{M}_{g,0}(\proj^1, d)$ be the moduli
stack of genus $g$, degree $d$ maps to $\proj^1$.
Consider the  $\com^*$ action on $\proj(V)=\proj^1$ as defined
in Section 2.
As before, there are canonical maps
$$\pi: U \rarr \overline{M}_{g,0}(\proj^1, d), 
\ \ \ \mu: U \rarr \proj(V)$$
where $U$ is the universal curve over the moduli stack.
Let $N$ denote the bundle $\oh_{\proj^1}(-1) \oplus 
\oh_{\proj^1}(-1)$.
Let
$$C(g,d)= \int _{[ \barr{M}_{g,0}(\proj^1,d)]^{\it{vir}}} 
c_{\rm{top}} (R^1 \pi_* \mu^* N).
$$
For each pair of linearizations $[\alpha, \alpha+1]$,
$[\beta, \beta+1]$, the virtual localization
formula yields an explicit computation of $C(g,d)$.

For general choices of linearization, $C(g,d)$
is expressed as a  complicated sum over connected 
graphs $\Gamma$ (see [GrP]) indexing the $\com^*$-fixed loci of 
$\overline{M}_{g,0}(\proj^1,d)$.
The vertices of these graphs lie over the
fixed points $p_1, p_2 \in \proj^1$ and are
labelled with genera (which sum over the graph to $g-h^1(\Gamma)$).
The edges of the graphs lie over $\proj^1$ and
are labelled with degrees (which sum over the 
graph to $d$). However, for the natural
linearization  $[0,1], [0,1]$,
a vanishing result holds: if a graph $\Gamma$
contains a vertex lying over $p_1$ of 
genus greater than 0 or
valence greater than
$1$, then the contribution of $\Gamma$ to $C(g,d)$
vanishes. As a result, the sum over graphs reduces to
a more manageable sum over partitions of $d$.
This linearization was found by Manin and used to
compute $C(0,d)= 1/d^3$ in [Ma]. In [GrP], the same
choice was used to compute $C(1,d)=1/12d$.

A dramatic improvement occurs if the
linearization $[0,1], [-1,0]$ is chosen.
In this case, a stronger vanishing holds:
if a graph $\Gamma$ contains
{\em any} vertex of valence greater than
$1$, then the contribution of $\Gamma$ to $C(g,d)$
vanishes. Hence, contributing graphs have exactly 1 edge.
The graph sum then reduces simply to
a sum over partitions $g_1+g_2=g$ of the genus.
The localization formula  yields the following result
for $g\geq 0$:
\begin{equation}
\label{rew}
C(g,d) = d^{2g-3} \sum_{\substack{g_1+g_2=g\\g_1,\,g_2 \geq 0}} 
b_{g_1} b_{g_2}
\end{equation}
where $b_g$ is defined by (\ref{qwq}).
In particular, the computations of $C(0,d)$ and $C(1,d)$
now
require no series manipulations of the type pursued
in [Ma], [GrP]. 
Note equation (\ref{rew})
implies
\begin{equation}
\label{rews}
\sum_{g\geq 0} C(g,1) t^{2g} = f_0(t)^2= f_1(t).
\end{equation}
In Section 4.2, the formula (for $g\geq 1$)
$$\sum_{\substack{g_1+g_2=g\\g_1,\,g_2 \geq 0}}
b_{g_1} b_{g_2} = \frac{|B_{2g}|}{2g}
\frac{1}{(2g-2)!}$$
will be proven from Theorem 2 and Bernoulli identities
to complete the proof of Theorem 3.

\subsection{Identities} Recall,
the Bernoulli numbers $B_m$ are defined by the
series expansion
\begin{equation}
\label{defff}
\beta(t)=\dfrac{t}{e^t-1} = \sum\limits_{m=0}^{\infty}B_m\dfrac{t^m}{m!}.
\end{equation}
We start by computing $b_g$ explicitly in terms of Bernoulli
numbers.

\begin{lm}
\label{bex}
$$
\frac{t/2}{\sin(t/2)}=1+\sum_{g\ge1}\frac{2^{2g-1}-1}{2^{2g-1}}
\frac{|B_{2g}|}{(2g)!}t^{2g}.
$$
\end{lm}

\bpf
This is well-known. We include a proof only for the reader's convenience.
\begin{eqnarray*}
\frac{t/2}{\sin(t/2)}&=&
\frac{it}{e^{it}-1}e^{it/2}=
\frac{it}{e^{it/2}-1}-\frac{it}{e^{it}-1}
=2\beta(it/2)-\beta(it)   \\
&=&2-\frac12it-\sum_{g\ge1}\frac1{2^{2g-1}}\frac{|B_{2g}|}{(2g)!}t^{2g}
-\left(1-\frac12it-\sum_{g\ge1}\frac{|B_{2g}|}{(2g)!}t^{2g}\right)  \\
&=&1+\sum_{g\ge1}\frac{2^{2g-1}-1}{2^{2g-1}}\frac{|B_{2g}|}{(2g)!}t^{2g}.
\end{eqnarray*}
\epf

By Theorem 2, we see (for $g\geq 1$)
$$b_g = \int_{\overline{M}_{g,1}} \psi_1^{2g-2}\lambda_g =
\frac{2^{2g-1}-1}{2^{2g-1}}
\frac{|B_{2g}|}{(2g)!}.$$
The series $f_1(t)= f_0^2(t)$ is determined by the
following lemma.

\begin{lm} \label{bbbq}
$$
\sum_{h=0}^{g}b_hb_{g-h}=\frac{|B_{2g}|}{2g}
\frac1{(2g-2)!}\,.
$$
\end{lm}

\bpf
Set $\beta_g=(2-2^{2g})\dfrac{B_{2g}}{(2g)!}$.
The identity to be proved (for $g\ge1$) is then
\begin{equation}
\label{edef}
2\beta_g+\sum_{h=1}^{g-1}\beta_h\beta_{g-h}=
-\frac{2^{2g}}{2g}\frac{B_{2g}}{(2g-2)!}\,.
\end{equation}
Since $\sum\limits_{g=0}^{\infty}\beta_gx^{2g-1}=\dfrac1{\sinh(x)}$
and $\sum\limits_{g=0}^{\infty}\dfrac{2^{2g}B_{2g}}{(2g)!}x^{2g-1}=\coth(x)$,
equation (\ref{edef}) is an immediate consequence of
$(\coth x)'= -\sinh^{-2}x. 
$
\epf

Lemma \ref{bbbq} yields the equality
$$f_1(t)= 1 + \sum_{g\geq 1} \frac{|B_{2g}|}{2g} \frac{t^{2g}}{(2g-2)!}.$$
This result together with  equations (\ref{rew}-\ref{rews}) completes
the proof of Theorem 3. \qed

\subsection{Proof of Theorem 4}

The equality (for $g\geq 2$)
$$
\int_{\overline{M}_g}\lambda_{g-1}^3=\frac{|B_{2g}|}{2g}
\frac{|B_{2g-2}|}{2g-2}\frac1{(2g-2)!}
$$
now  may be established by manipulating Mumford's
Grothendieck-Riemann-Roch formulas and using  Lemma \ref{bbbq}.

\vspace{+10pt}
\bpf The formula 
$\sum_{k\ge1}(-1)^{k-1}(k-1)!\,\ch_k(V)t^k
=\log(\sum_{k\ge0}c_k(V)t^k)$ gives for $V=\hodge$
$$
\sum_{k\ge1}(-1)^{k-1}k!\,\ch_k(\hodge)t^{k-1}=
\left(\sum_{k=1}^g k\lambda_kt^{k-1}\right)
\left(\sum_{k=0}^g \lambda_k(-t)^k\right)
$$
since $c(\hodge)^{-1}=c(\hodge^{*})$.
(Note that both sides are even polynomials in $t$.)
In particular $(2g-3)!\,\ch_{2g-3}(\hodge)=
(-1)^{g-1}(3\lambda_g\lambda_{g-3}-\lambda_{g-1}\lambda_{g-2})$
so that 
$$
\lambda_g\lambda_{g-1}\lambda_{g-2}=
(-1)^g(2g-3)!\,\lambda_g \ch_{2g-3}(\hodge).
$$
Mumford's formula [Mu] for $\ch(\hodge)$ gives
$$
(2g-3)!\,\ch_{2g-3}(\hodge)=
\frac{ B_{2g-2}}{2g-2}\left[
\kappa_{2g-3}+\frac12\sum_{h=0}^{g-1}i_{h,*}
\left(\sum_{i=0}^{2g-4}\psi_1^i(-\psi_2)^{2g-4-i}\right)\right].
$$
Since $\lambda_g=0$ on $\Delta_0$ 
while $i_h^*\lambda_g ={\text{pr}}_1^*\lambda_h{\text{pr}}_2^*\lambda_{g-h}$
for $h>0$, this implies
$$
\int_{\overline{M}_g} \lambda_{g-1}^3=
\int_{\overline{M}_g} 2\lambda_g\lambda_{g-1}\lambda_{g-2}=
\frac{|B_{2g-2}|}{2g-2}
\left[2b_g+\sum_{h=1}^{g-1}b_hb_{g-h}\right]
$$
(where the first equality follows from $c(\hodge) c(\hodge^*)=1$).
Hence, it remains to prove 
$$
\sum_{h=0}^{g}b_hb_{g-h}=\frac{|B_{2g}|}{2g}
\frac1{(2g-2)!}.
$$
But, this is Lemma \ref{bbbq}. \epf

\subsection{The Virasoro prediction for $c_g$}
We include here D. Zagier's proof 
of the prediction  (for $g\geq 1$):
$$\Bigl(
\sum_{k=1}^{2g-1} \frac{1}{k} \Bigr) b_g = c_g + \frac12 \sum_{g=g_1+g_2,
g_i>0}
\frac{(2g_1-1)! (2g_2-1)!}{(2g-1)!} b_{g_1} b_{g_2}.$$
{From} Theorem 2, we obtain
$$\sum_{g\geq 1} c_g t^{2g}= 
\Big( \frac{t/2}{\sin(t/2)} \Big) 
\cdot \log\Big( \frac{t/2}{\sin(t/2)} \Big).$$
Lemma 3 below (together with Lemma \ref{bex}) expresses
$c_g$ in terms of Bernoulli numbers.
Then, the Virasoro prediction for $c_g$ is equivalent to
an identity among Bernoulli numbers proven in Lemma 4.
 
\begin{lm} 
$$
\log \left( \frac{t/2}{\sin(t/2)} \right) =
\sum_{k\ge1} \frac{|B_{2k}|}{(2k)(2k)!}t^{2k} .
$$
\end{lm}

\bpf
Let $f(t)=\dfrac{t/2}{\sin(t/2)}$. It suffices to prove
\begin{equation}
\label{rrrt}
t\frac{f'(t)}{f(t)}=\sum_{k\ge1}\frac{|B_{2k}|}{(2k)!}t^{2k}.
\end{equation}
By definition (\ref{defff}), the right  side of (\ref{rrrt}) equals 
$1-\dfrac12it-
\dfrac{it}{e^{it}-1}.$
The left side equals $1-\dfrac{t}2\cot(t/2)=1-i\dfrac{t}2
\dfrac{e^{it}+1}{e^{it}-1}$. \epf

\begin{lm} \label{Zag}
\begin{eqnarray*}
\lefteqn{
\left(\sum_{l=1}^{2g-1}\frac1l\right)\frac{2^{2g-1}-1}{2^{2g-1}}
\frac{|B_{2g}|}{(2g)!}=
\sum_{k=0}^{g-1}
\frac{|2^{2k-1}-1|}{2^{2k-1}}
\frac{|B_{2k}|}{(2k)!}
\frac{|B_{2g-2k}|}{(2g-2k)(2g-2k)!} } \cr
&&\qquad\qquad\mbox{}+\frac12\sum_{\substack{g_1+g_2=g\\g_1,\,g_2>0}}
\frac1{(2g-1)!}
\frac{2^{2g_1-1}-1}{2^{2g_1-1}}
\frac{2^{2g_2-1}-1}{2^{2g_2-1}}
\frac{|B_{2g_1}|}{2g_1}
\frac{|B_{2g_2}|}{2g_2}.
\end{eqnarray*}
\end{lm}

\noindent {{\em Proof} (Zagier).}
\def\b{\beta}\def\s{\sinh}\def\f{\frac}\noindent
Set $\b_g=(2-2^{2g})\,\dfrac{B_{2g}}{(2g)!}\,$. The identity to be proved
is $a(g)+b(g)=c(g)$, where
\begin{eqnarray*}
&&a(g):=\biggl(1+\f12+\cdots+\f1{2g-1}\biggr)\,\b_g\,,\\
&&b(g):=\sum_{n=1}^g\f{2^{2n}B_{2n}}{2n\,(2n)!}\,\b_{g-n}\,,\\
&&c(g):=\f12\sum_{\substack{g_1+g_2=g\\g_1,\,g_2>0}}
 \f{(2g_1-1)!\,(2g_2-1)!}{(2g-1)!}\,\b_{g_1}\b_{g_2}\,.
\end{eqnarray*}
Using the generating function identity
$\sum\limits_{g=0}^\infty\b_g\,x^{2g-1}=\dfrac1{\s x}$, we find
\begin{eqnarray*}
A(x)&:=&\sum_{g=1}^\infty\,a(g)\,x^{2g-1}
  =\sum_{g=1}^\infty\b_g\int_0^x\f{x^{2g-1}-t^{2g-1}}{x-t}\,dt  \\
 &=&\int_0^x\biggl[\f1{x-t}\,\biggl(\f1{\s x}-\f1{\s t}\biggr)
   +\f1{xt}\biggr]\,dt\,,\\
B(x)&:=&\sum_{g=1}^\infty\,b(g)\,x^{2g-1}
  =\f1{\s x}\,\sum_{n=1}^\infty\f{2^{2n}B_{2n}}{2n\,(2n)!}\,x^{2n}
  =\f1{\s x}\,\log\bigl(\f{\s x}x\bigr)\,,\\
\end{eqnarray*}
\begin{eqnarray*}
C(x)&:=&\sum_{g=1}^\infty\,c(g)\,x^{2g-1}
  =\f12\,\sum_{\substack{g_1+g_2=g\\g_1,\,g_2>0}}\b_{g_1}\b_{g_2}\,
   \int_0^xt^{2g_1-1}\,(x-t)^{2g_2-1}\,dt\\
  &=&\f12\int_0^x\biggl(\f1t-\f1{\s t}\biggr)
  \biggl(\f1{x-t}-\f1{\s(x-t)}\biggr)\,dt
\end{eqnarray*}
and hence
\begin{eqnarray*}
\lefteqn{2\,C(x)-2\,A(x)\;=\int_0^x\biggl\{
\biggl(\f1t-\f1{\s t}\biggr)\biggl(\f1{x-t}-\f1{\s(x-t)}\biggr)}\\
&&\mbox{}-\biggl(\f1{x-t}+\f1t\biggr)\biggl(\f1{\s x}+\f1x\biggr)
  +\f1{x-t}\,\f1{\s t}+\f1t\,\f1{\s(x-t)}\biggr\}\,dt\\
&&=\int_0^x\biggl(\f1{\s(t)\s(x-t)}-\f x{\s x}\,\f1{t\,(x-t)}\biggr)\,dt\\
&&=\f1{\s x}\,\log\biggl(\f{\s t}t\cdot\f{x-t}{\s(x-t)}\biggr)
 \biggl.\biggr|_0^x\;=\;2\,B(x)\,.
\end{eqnarray*}

\vspace{+10 pt}
\noindent Department of Mathematics \\
\noindent Oklahoma State University \\
\noindent Stillwater, OK 74078 \\
\noindent cffaber@math.okstate.edu 

\vspace{+5pt}
\noindent Institutionen f\"or Matematik \\
\noindent KTH \\
\noindent 100 44 Stockholm, Sweden \\
\noindent carel@math.kth.se 

\vspace{+10 pt}
\noindent
Department of Mathematics \\
\noindent California Institute of Technology \\
\noindent Pasadena, CA 91125 \\
\noindent rahulp@cco.caltech.edu
\end{document}